# Gender Gaps in the Mathematical Sciences:
# The Creativity Factor

by Theodore P. Hill and Erika Rogers

The underrepresentation of women in the mathematically-intensive sciences (hereafter, for brevity, referred to as the *hard sciences*) has been a concern in the United States for over half a century. After the women's movement began in the 1960s, gender gaps in many professional fields decreased, and some even reversed dramatically. Currently women comprise about half the M.D's, two-thirds of psychology Ph.D's, and three-quarters of veterinary medicine doctorates, more than seven times as high as in the 60s [CW3, p. 5]. In the hard sciences, however, the large gender gap favoring men has stubbornly persisted (see Figure 1), and many efforts have been made to determine why.

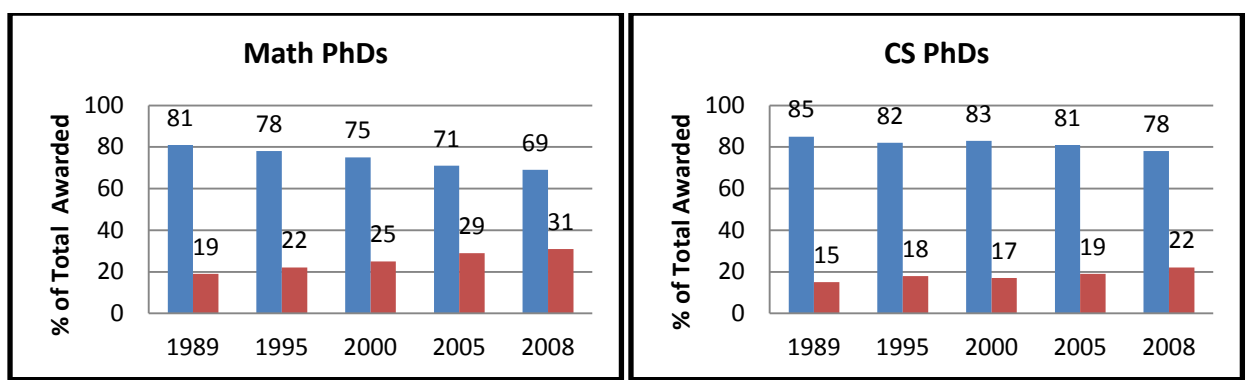

**Figure 1.** (Source: [NSF2]. Red shows percent women and blue shows percent men.)

Hundreds of millions of dollars of public funds are being devoted to understanding this particular gender gap "problem". For example, in addition to its many other programs that indirectly support women in science [NSF2], the goal of the U.S. National Science Foundation's (NSF) special ADVANCE program (Advancement of Women in Academic Science and Engineering Careers) is "to increase the representation and advancement of women in academic science and engineering careers" [NSF1]. In the past ten years the ADVANCE program alone has awarded over $130 million of public funds to this cause.

The NIH (National Institutes of Health) also supports research in the scientific gender gap field, such as a current grant of $1.4 million to two faculty researchers for a single three-year study entitled *Assessing and Reducing Gender Bias in STEM* [Science, Technology, Engineering and Mathematics]. Many other government agencies and state universities also contribute significant public resources to addressing this gender gap in science, and continue to solicit further applications for such awards (e.g., [CIWS, NSF3, UCB]). During her recent presidential campaign, Hillary Clinton argued that "women comprise 43 percent of the workforce, but only 23 percent of scientists and engineers", urging the government to take "diversity into account when awarding education and research grants" [CW3, p 54]. However, despite these efforts, it appears that only a little progress has been made, and the causes still elude us.

**Prevailing Theories**

The recent flurry of books, anthologies, survey articles, and book reviews on the gender gap in science [e.g., Col, CW1-4, CWB, GK, KM, Sc, Sp] includes a new three-year study summarizing the findings of over 400 research articles and "approximately 20 meta-analyses (and several meta-analyses of meta-analyses)" [CWB, p 219]. After careful analysis of this huge body of literature, Cornell developmental psychologists Stephen Ceci and Wendy Williams rule out discrimination as a significant factor, even describing evidence of reverse-discrimination. As confirmed by many of our own colleagues of both genders, there simply are no queues of girls and women striving to enter careers in mathematics, computer science and other hard sciences, and being turned away. If anything, just the opposite is true.

The new study [CW3] reviews the standard gender gap arguments such as pipeline issues, motherhood, the "people" versus "things" explanation, and the "greater male variability hypothesis" (that men and women are of equal average ability, but that the variance of men is higher – hence more idiots and more geniuses). Then, the authors find that the evidence is not consistent with the gender gap being largely a consequence of biological sex differences [CW3, p 180], and that social factors are also not compelling [CW3, p 183]. The meta-analysis concludes:

> we believe that the entire corpus of research reduces to a single large effect coupled with a host of smaller effects. The largest effect concerns *women's choices and preferences* – their preference for non-math careers over careers in engineering, physics, mathematics, operations research, computer science and chemistry [CW3, p 179-80, emphasis added].

That conclusion is neither new, nor widely accepted among scientists themselves. For example, in her introductory comments to the 1999 landmark report by the Committees on Women Faculty in the School of Science at MIT, chair of faculty Lotte Bailyn wrote "Our first instinct is to deny that a problem exists (if it existed, it would surely have been solved by now) or to blame it on the pipeline or the circumstances and choices of individual women" [Ba]. And as Oberlin mathematician Susan Colley opined, "I find the issue of "personal choice" to be more complicated than Ceci and Williams seem to" [Col].

After reaching the "women's preference" conclusion, the authors of *The Mathematics of Sex* "challenge those with different views to present evidence they believe we ignored or misinterpreted" [CW3, p 15]. The studies [CW1-4, GK] include input from scores of sociologists, psychologists, educators, government appointees, biologists, boys and girls and male and female rats, but apparently few, if any, actual hard scientists. Indeed, Ceci and Williams admit that "We do not know what it takes to be a successful math, engineering, or physics professor, or a chemist or computer scientist" [CW3, p. 83]. It is one of the goals of this article to provide evidence that some hard scientists feel was ignored, and to suggest a theory-driven partial explanation for the gender gap in the hard sciences.

Even if the "women's preference" conclusion is accepted, the original question of "Why?" remains unanswered, and, perhaps more importantly, so does the question of what could or even *should* be done about it. Do the majority of women *prefer* not to go into the hard sciences because of their own limitations in either aptitude or attitude (i.e., they simply don't have the talent, or they think they don't have the talent), or because there's something intrinsically unappealing to them about these fields? And what about the women who *do* go into these fields, and then leave? The issue of raising children simply does not account for the smaller influxes and larger exoduses observed in hard science careers over others. Is there some other important common factor that should be considered?

**The Creativity Factor**

The self-described top researchers in the gender gap in science [CW1] seem to have completely ignored an important and compelling factor. In spite of acknowledging up front "the kind of intense, *highly creative thinking required of mathematicians*" ([CW3, p x], emphasis added), they have omitted the well-studied issue of gender differences in creativity. In ignoring the creativity factor, the science gender gap experts have greatly underestimated the potential importance of a completely different set of both biological and societal factors which may "conspire to limit talented women and girls" [CW3]. Consequently, decision makers are thereby missing significant opportunities for constructive improvements.

If the topic of gender differences in the hard sciences "has initiated such strong and impassioned reactions that it has not always been a suitable topic for dinner conversations" [Col, p 379], the topic of gender differences in creativity is downright inflammatory. In mathematician Reuben Hersh's words about another mathematics overrepresentation issue, "Too ticklish, too much chance to be misunderstood, or give offense, or get in trouble one way or the other" [He]. Creativity experts John Baer and James Kaufman freely concede that the gender difference topic is "a difficult arena in which to conduct research" [BK, p 75].

The notion of creativity itself is difficult, and meta-analyses of the field such as [BK, P2] do not even attempt to provide a clear definition. Among the scores of characterizations in the psychology literature, there is no single, authoritative definition of creativity [FBAM]. There are also many different measures of creativity, such as the Torrance Tests of Creative Thinking and various related tests, but "For at least 25 years a debate has raged over the validity of these tests as measures of creativity" [BK, p 79].

In industry and academia, creativity is also undefined, but is often measured in terms of quantity and quality of various outcomes: patents, numbers of papers or citations, grants received, profitable innovations, and prizes. In the special case of mathematics, it is widely accepted that the highest prize is the Fields Medal, officially known as the "International Medal for Outstanding Discoveries in Mathematics", and this may perhaps be viewed as reflecting the very highest level of creativity in the subject. But as musician Matt Callahan said, "Something as porous as creativity defies definition, resists quantification and refuses access to those who seek to possess it like a Thing" (as cited in [P2, p 6]).

However it is defined, the subject of creativity has a long history of research. While it is often associated with art and music, creativity is clearly also a key factor in high performance mathematics and hard sciences [Ad], and therefore merits serious attention as a contributing element to success in these fields. In the literature on gender gaps in science, however, the words "creative" or "creativity" do not even appear in the indices of either [CW3] or [GK].

Creativity also does not appear among the plethora of studies reviewed in [CW3] of other differences between the genders – differences in abstract thinking and reasoning, in academic scores, in brain size and structure, in early activities and math competence, in journal article writing, in reasoning abilities, in social skill development, in verbal processing, and, above all, in spatial reasoning. As Stanford historian Londa Schiebinger wrote, "We as a society do support endless studies of sex differences in spatial perception" [Sc]. The role of gender gaps in creativity, however, is essentially missing in the studies of gender gaps in science.

In spite of the fact that neither the topic of creativity nor divergent thinking appears in Halpern's *Sex Differences in Cognitive Abilities* [BK, p 76], gender differences in creativity are well-studied (e.g., the survey [BK] contains over 180 references) and are widely accepted [BK, P1, P2, RB]. Two important facets of creativity are sometimes distinguished, namely *creative ability* and *creative achievement*. Findings from studies on gender differences in individual creative ability, including standard creativity tests, self-reports, personality tests, and teacher/peer assessments, are all over the map. The results are contradictory and inconclusive. Interestingly, the list of studies concluding that girls and women are *more* creative than boys and men are more numerous than those with the opposite conclusion [BK, pp 80-87].

On the other hand*,* there seems to be broad consensus that there are gender differences in creative achievement "at the highest levels, as judged by the experts in their respective domains, with men dominating most fields" [BK, p 97], and that women "appear more interested in the creative process itself than in its end-product" [RB, p 100-101]. As Ashland University creativity expert Jane Piirto puts it,

> The women's movement began in the 1960's…. Why have we not begun to see a more equal ratio of successful women to men in creative fields? Where are the publicly and professionally successful women visual artists, musicians, *mathematicians*, *scientists*, composers, film directors, playwrights, and architects?...It seems that the only creative fields where women are equally known as men are creative writing and acting…. [P1, p 142, emphasis added].

For example, since the Fields Medal was inaugurated more than 75 years ago, 52 awards have been made, and not a single one has been to a woman. With respect to the more general scientific community, Ceci and Williams observed that the overrepresentation of men is larger in disciplines requiring more mathematics, and raised the question why this problem should be so much worse for math-intensive fields than other high-powered professions [CW3, p. 104].

Since both society and experts seem to agree that there is a great difference between women and men in creative achievement at the highest levels, and since gender gap experts also acknowledge that mathematics requires highly creative thinking, it is reasonable to wonder whether a significant factor in explaining the dearth of women in the hard sciences may also have to do with gender differences in creative achievement.

Nearly as many girls now take mathematics in college as boys, and girls get better grades in mathematics [CW3, p 30]. The big drop off in numbers for the hard sciences apparently begins at the point of graduate work at the Ph.D. level, and continues through the tenure-track procedure and out the other end where even successful women scientists are voluntarily leaving these fields in larger numbers than men [CW3, p 7]. Why should this be? One consideration is that the Ph.D. dissertation is where the most creative, original and challenging academic work *begins*.

Piirto, for example, cites several studies indicating that differences between creatively gifted males and females seem to come "in the *choices that they make after college*, a time when commitment and regular effort in the field of creativity matters" [P1, emphasis added]. She also observed that girls do not show less creative achievement until *after high school and college* [as cited in BK, p 94, emphasis added]. Thus, recent conclusions about gender differences in mathematics performance among boys and girls, such as the studies of hundreds of thousands of fourth and eighth graders from forty-eight countries reported in [KM], seem to miss the critical period when creative achievement begins to kick in. Are so many post-college women choosing not to proceed in the hard sciences for the same reasons that they seem to shy away from careers in other highly creative fields?

**Factors Related to Creativity**

Creative achievement is seen to be enhanced by many factors, and there is evidence of gender differences in several of these (e.g., [P1]). We take a brief look at three of these contributing factors. First, men are often seen to be more playful than women, more immature than women [Cr], and in general happier than women, as two recent studies reported in the *New York Times* found [Le]. If we look at Gary Larson's cartoons of scientists "at work", such as the famous one with two balding male scientists in lab coats, one working on an atomic bomb while his buddy is sneaking up behind him about to burst a paper bag of air, the humor is immediate. Would this seem as funny if the scientists were two middle-aged, nerdy women instead of men?

Yet play has been recognized as an important catalyst for the creative mind, not only for children, but also for adults working in organizational settings [MR]. A colleague's anecdote from industry illustrates this point. When she was a math intern at Bell Labs, she was shocked at how many famous mathematicians just sat around playing cards and Go all day. Then suddenly one day, a player would drop his hand of cards on the table and excitedly start talking science and drawing diagrams in the air. The others soon chimed in, and the idea they spawned during the next hour led to an invention that paid all their annual salaries. Next day, more fun and games. It is easy to believe this scenario, but hard to imagine a group of highly educated and creative women acting the same way.

Another factor associated with success in science, perhaps especially in the "laboratories in the mind" of hard sciences, is curiosity. Here too, according to some studies, men sometimes appear to have an advantage. For example, George Mason psychologist Todd Kashdan and colleagues conducted four studies based on their Curiosity and Exploration Inventory (CEI), which comprises two dimensions: "exploration (appetitive strivings for novelty and challenge) and absorption (full engagement in specific activities)". Men reported statistically significant greater CEI exploration scores ($p < .01$), greater absorption scores ($p < .05$) and greater CEI-total scores ($p < .01$) on one of the samples, and similar trends were found for three other samples, although none was significant [KRF, p 295]. Since "curiosity functions as an adaptive motivational process related to the *pursuit of novelty* or challenge" [GL, p 236, emphasis added], this too might contribute to understanding gender differences in creative achievement.

A third factor positively associated with creativity (e.g., [RH]), and consequently with success in the hard sciences, is the willingness to take risks, and to accept rejection and failure. When a writer asked Thomas Watson, founder of IBM, for the secret to success, Watson's famous answer was "Double your failure rate". Today is no different than in Galileo's era in that successful scientists routinely experience rejection – rejection of papers for publication, of positions at top universities, of Nobel prizes and Fields Medals and scores of lesser accolades. The ability to persevere in the face of repeated rejection applies to a wide variety of disciplines, but the humor in Sidney Harris's cartoon showing a bearded professor hunched over his desk in an office labeled DEPT OF MATHEMATICS & FRUSTRATION is evident, even to a lay person. Would a sign saying Department of Biology and Frustration be seen as equally humorous, except perhaps by biologists? The mathematically-intensive sciences are particularly brutal with respect to what constitutes a successful result: "Mathematics, like chess, requires too direct and personal a confrontation to allow graceful defeat" [Ad, p 3].

Men are viewed as better able to accept rejection, and in that respect, the mathematically creative personality is akin to that of his artist colleagues. Former chair of Washington State's Department of Dance and Theater, Laurilyn Harris, noted

> The profession of artist *demands an extraordinary commitment in terms of willingness to take rejection*, to live in poverty, and to be field independent. Those are traits of committed males, but not of committed females, who usually choose careers as art educators, but not as artists (as cited in [P1], emphasis added).

While numerous studies on gender differences in risk taking (or risk aversion) "support the idea that male participants are more likely to take risks than female participants" [BMS, p 377], many questions remain unanswered. For example, is there a relationship between the kinds of artistic career choices mentioned above and those of women preferring not to go into the hard sciences? Are such choices impacted by risk-taking characteristics, and if so, can (or should) they be modified?

In all of these examples, we see that there are common elements across the creative fields, whether artistic, scientific, or mathematical. Thus, in order to avoid tunnel vision, it may be useful for researchers as well as decision makers to weigh possible solutions to the gender gap problem in the hard sciences against the broader context of the creative spectrum. For example, although many agree that intense and focused commitment in general is a key to both artistic and scientific productivity, "[s]ome of the proponents of gender equity [have made] demands to *abolish the obsessive and compulsive work ethic of successful scientists* that universities reward" [CW3, p 195, emphasis added]. Would they also abolish the work ethic of successful artists, writers, composers, and chefs?

## Nature or Nurture?

If men are more creatively productive at the top end, including the hard sciences, why might that be? One controversial argument supporting significant biological gender differences in creativity is simple. Experts estimate that humans, as a species, have been hunter-gatherers for all but 600 of their 10,000 generation history [Li]. Some studies suggest that "male specializations in hunting and making artifacts may have been more cognitively demanding than female specializations in gathering and child rearing" [CWB, p 237]. According to this argument, over a period of 9000 generations, evolution could also have contributed to a gender gap in creative thinking (as well as in more obvious traits such as size, aggressiveness, etc.). However, other researchers, including Ceci *et al*, have found that "the available evidence is insufficient to determine the impact of evolution on sex differences in cognitive ability, although it presents intriguing suggestions" [CWB, p 237].

What we do know is that male and female human brains are now physically different (e.g., [BL, CW3]), and these physical differences may also be reflected in different thinking processes. However, since biology cannot readily be changed, whether or not the gender gap in creative achievement is innate has limited use, except perhaps to suggest flexibility in what gender "equity" means numerically. No one seems to argue for exactly fifty-fifty, but 70-30 is seen as a problem in some scientific professions, but not in others such as scientific medical research (see Figure 2).

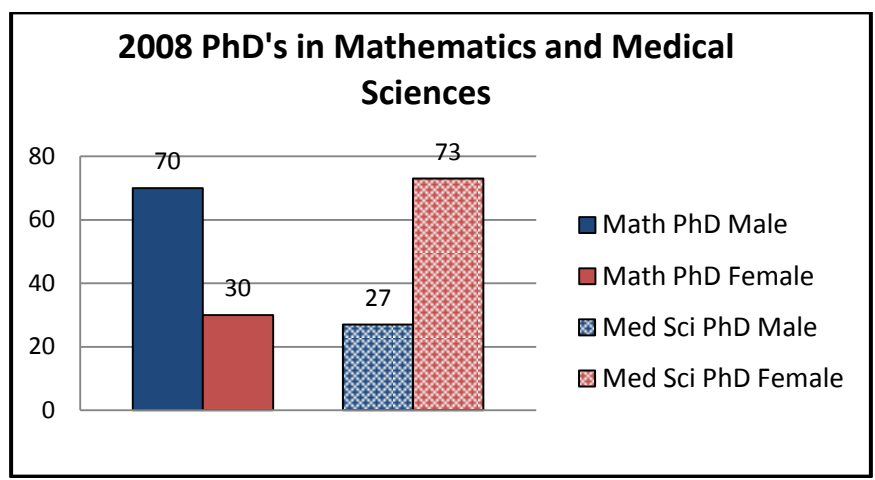

**Figure 2.** (Source: [NSF2])

Society, on the other hand, *can* be changed, and thus it is far more important to recognize that there are significant cultural and societal reasons for the gender gap in creative achievement. Girls are often steered away from "unladylike" playful behavior, from getting dirty to tearing devices apart. As Syracuse University developmental psychologist Alice Honig found, "boys are frequently allowed more freedom of movement, more permission to cross streets and roam further in neighborhoods, more indulgence for climbing and jumping" [Ho, p 115].

Equally important, girls do not have the same number of creative role models - composers, architects, scientists, chefs, inventors, playwrights and film directors - as boys. Even in today's computerized society, who are the hackers, and the inventors of Google, Facebook, and the Internet itself? There are certainly more female role models in mathematics now than in the 1960s, at least at the level of professors and researchers, but at the very highest levels, such as the Fields Medal or Nobel prizes in the hard sciences, the situation is as bad as ever.

On the other hand, girls and women are still more heavily burdened with family responsibilities and expectations that may compete with their choices to lead creative lives. Gender differences in terms of men's participation in domestic labor and child care, as well as conflicting societal approval of caregiver roles add to this mix. Solutions such as the provision of childcare facilities in the workplace (e.g., [Sp]) will certainly benefit women and men across all disciplines. In addition, however, there are some steps that can be taken to specifically support those engaged in highly creative fields.

## Constructive Opportunities

A major goal of society, presumably, is to determine how best to utilize the talents of its individuals for the greater good of that society. It has also been suggested that diversity contributes to a richer mix of ideas, inventions, innovations, and problem solutions. Simply addressing gender differences in creativity is certainly not a panacea for the gender gap in the hard sciences. However, recognizing that intensive creative thinking and achievement is an important component of success in these disciplines opens the door to a wealth of opportunities.

Baer and Kaufman argue that there is at least one over-arching reason why women's creative achievement has lagged in almost all fields, and that is the relative lack of environments conducive to developing expertise [BK, p 77]. But environments are one thing that can readily be improved, sometimes through relatively simple means. For instance, as mentioned earlier, simple play is frequently the catalyst to new ideas. This is aptly demonstrated in Tim Brown's lively "Serious Play" TED lecture using Finger Blasters [Br]. It is even more strongly supported in the discussion of play and creativity in the workplace by Charalampos Mainemelis and Sarah Ronson of the London Business School: "when play is woven into the deep fabric of organizational life it can transform the very nature of their products and work processes" [MR, p 83-84].

Thus one constructive idea for enhancing the creative output of an organization is simply to install playrooms in the workplace such as those at IDEO and Google. A room with computer and board games, Legos and modeling clay, whiteboards and colored pens might well lead to unexpected dividends in discoveries in the hard sciences, by both sexes. This approach may also contribute to a general climate of happiness at work, which some say is the number one productivity booster.

Of course, the current budget cuts in colleges and universities make new expenditures difficult, but they could begin, for example, with changes in faculty coffee rooms. Of the scores of mathematics and computer science faculty lounges we have visited here and abroad, every single one could be transformed into a much more fun place with the addition of a few carefully selected games and toys. Whether immersion for several hours, or days, or weeks in Google-style playrooms generates creativity, or releases it, seems immaterial if the bottom line is more creative output.

A more direct approach to increasing creative output is exemplified by the Hasso Plattner Institute of Design at Stanford. Recognizing that it is much harder to teach creativity than core subjects, the *d-school*, as it is called, provides an "innovation hothouse" where students engage in "ideation" [An]. The goal of the *d*-school, in short, is to teach imagination. It has already spawned impressive tangible achievements such as inexpensive, solar-powered lamps for the rural poor in the developing world, and graduate students compete fiercely to get into the program. Thus if an institution's goal is to increase the success of hard scientists of a particular gender (or race, say), it could set up its own innovation hothouse, and especially encourage students in those groups to attend. Again, as with departmental "playrooms", it could be relatively cheap and easy to devote one classroom to an innovation hothouse, where graduate students and faculty could spend off hours when they need inspiration.

Related to the innovation labs is the idea for institutions to provide some form of explicit "failure training" for budding hard scientists who are struggling with the many rejections and dead-end ideas that are especially common in these fields. A key goal of failure training is to teach people to step back from the disappointment of an idea or investigation that didn't pan out, accept that the expenditure of time and effort was still worth it, and try to determine what was learned or valuable that can be taken away from this experience. The important thing is to keep at it, even after temporary setbacks.

Since women are perceived as less willing to risk rejection, such workshops could especially benefit women scientists. Stanford's *d*-school, for example, immerses students in what they call a "constant churning of rethinking, repurposing and recommitting, even when they've been battered by a series of early failures…if someone's creative energy gets drained, there's a shoes-off white room to retreat to, where scrawling on the floor and walls may stir a breakthrough" [An]. Similarly, the renowned Isaac Newton Institute for the Mathematical Sciences in Cambridge, England, even has blackboards in the restrooms.

Finally, we want to report a very interesting and elegant idea we learned during a recent visit to the mathematics department at the United States Military Academy at West Point. The female math professors at USMA, recognizing that mathematics research often requires intense solitary thinking and concentration, set aside times to meet together in the library, where they occupy a room and sit down together, each working silently on her own research. Although this is not a "playful" activity, it reveals another dimension of how to shape an environment conducive to creative output. This program, which they dubbed GDR for "Girls Do Research", has been such a success that the USMA male math professors, in an attempt to play catch-up, established GDR2 (Guys Do Research Too).

Many of these ideas – playrooms, innovation labs, failure training, GDR groups – can also benefit other disciplines, of course, but they might prove especially effective for stimulating research output and success in areas requiring intense analytical creativity. Bell labs, Google and Stanford seem to think so.

**Conclusions**

In America today, gender gaps exist in the sciences. Some, such as mathematics, have persisted over generations, while others such as biology have reversed. While the continued underrepresentation of women in the hard sciences is viewed by many as a problem worth devoting substantial public resources to solve, the overrepresentation of women in the medical sciences apparently is not.

The hundreds of studies of gender gaps in science have ignored what we feel to be an important factor, namely gender gaps in creativity and related traits such as playfulness, curiosity, commitment and willingness to take risks. If our society chooses to continue to expend significant resources for psychological studies of gender gaps in science, at least some should be aimed at understanding essential creativity factors.

In the meantime, we feel that changes enhancing and encouraging a "culture of creative opportunity" for students *and* faculty could be implemented effectively and quickly within current academic environments, particularly those with a view to improving women's representation in the hard sciences. Perhaps the funding agencies could consider spending the next $130M by giving one million dollar grants to fund *d-schools* or Google-like playrooms at 130 institutions. This may not directly increase the relative creative achievement of women in the hard sciences, of course, but it seems worth a try, and if the result is better science, it will serve society nonetheless.

Are there any quick fixes to the gender gap "problem" in the hard sciences? One possible solution, of course, is to try to lure talented women from other fields into the hard sciences. But even if scientific talents *were* transferable from one field to another, Cornell psychologist Susan Barnett asked, "is it more valuable to encourage women to shift from their dominance in fields of biology to mathematics, so they can end up working on a search algorithm for Google rather than on a cure for AIDS?" [CW3, p 57]. On the other hand, a creative environment demonstrating how a Google search algorithm could facilitate a cure for AIDS might convince women who do show talent in several areas to opt for a hard sciences career.

A final idea is to jumpstart outstanding women who are already established hard scientists. As one of our colleagues suggested with a wide grin, maybe ADVANCE could identify 130 highly creative American women hard scientists, and give each a million dollars to prove theorems, invent technology or solve engineering problems. Then just sit back and watch the sparks fly!

**Acknowledgements.** The authors are grateful for excellent suggestions and comments from many friends and colleagues, the reviewers, and especially editor Marjorie Senechal. We also thank Charles Day for his *Physics Today* blog and Stephen Ceci for his comments.

## About the Authors

The first author is professor emeritus of mathematics at the Georgia Institute of Technology, currently Research Scholar in Residence at the California Polytechnic State University in San Luis Obispo. He holds a B.S in engineering from West Point, an M.S. in operations research from Stanford and a PhD in mathematics from Berkeley. His main research interest is in the mathematical theory of probability, especially optimal stopping, fair division problems, and Benford's law.  He served on two NSF panels in Washington for women and minorities in mathematics.

The second author is retired from the California Polytechnic State University, where she was a full professor of computer science, director of the university Honors Program, and faculty for the Adult Degree Program. She holds a B. Math. from the University of Waterloo, a PhD in computer science from the Georgia Institute of Technology, and an M.A. in Depth Psychology from Pacifica Graduate Institute. Her research interests include human-centered computing and usability, human-robot interaction, and human-information interaction. She currently works as a consultant and academic research coach.